\UseRawInputEncoding
\documentclass[12pt, oneside]{amsart}
\usepackage{geometry}
\numberwithin{equation}{section}
\usepackage{graphicx}
\usepackage{amsmath,amssymb,amsthm,amscd}
\usepackage{ascmac}
\usepackage{indentfirst}
\usepackage{cases}
\usepackage[all]{xy}
\usepackage[OT2,T1]{fontenc}
\DeclareMathOperator{\Gal}{\mathrm{Gal}}
\newtheorem{thm}{\textbf{Theorem}}[section]
\newtheorem*{thm*}{Theorem}

\newtheorem*{conj*}{\textbf{Conjecture}}
\newtheorem{defn}[thm]{\textbf{Definition}}
\newtheorem{prop}[thm]{\textbf{Proposition}}
\newtheorem{lem}[thm]{\textbf{Lemma}}
\newtheorem{sublem}[thm]{\textbf{Sublemma}}

\newtheorem{rem}[thm]{Remark}

\newcommand{\ilim}[1][]{\mathop{\varprojlim}\limits_{#1}}

\newcommand{\e}[1][]{\varepsilon}
\newcommand{\g}[1][]{\gamma_n}
\newcommand{\ka}[1][]{\kappa}
\newcommand{\G}[1][]{\Gamma}
\newcommand{\Gn}[1][]{\Gamma_n}
\newcommand{\gh}[1][]{\widehat{\gamma_n}}
\newcommand{\Gkn}[1][]{G_{K_n}}
\newcommand{\f}[1][]{\varphi}
\newcommand{\A}[1][]{A_{\mathrm{crys}}}
\newcommand{\Ap}[1][]{\widetilde{A}^{+}}
\newcommand{\Ac}[1][]{\widetilde{A}}
\newcommand{\Ak}[1][]{A_{K_n}}
\newcommand{\B}[1][]{B^{+}_{\mathrm{crys}}}
\newcommand{\Bp}[1][]{\widetilde{B}^{+}}
\newcommand{\Bc}[1][]{\widetilde{B}}
\newcommand{\Bk}[1][]{B_{K_n}}
\newcommand{\Q}[1][]{\mathbb{Q}}
\newcommand{\Z}[1][]{\mathbb{Z}}
\newcommand{\C}[1][]{\mathbb{C}}
\newcommand{\ch}[1][]{\chi_{\cyc}}
\newcommand{\cyc}[1][]{\mathrm{cyc}}
\newcommand{\F}[1][]{\mathbb{F}_p}
\newcommand{\z}[1][]{\zeta_{p^n}}
\newcommand{\Ep}[1][]{\widetilde{E}^{+}}
\newcommand{\E}[1][]{\widetilde{E}}
\newcommand{\Ok}[1][]{\mathcal{O}_{K}}
\newcommand{\Res}[1][]{\mathrm{Res}}
\newcommand{\Log}[1][]{\mathfrak{L}}
\newcommand{\Tr}[1][]{\mathrm{Tr}}
\newcommand{\tens}[1]{%
\mathbin{\mathop{\otimes}\limits_{#1}}%
}
\title[]{On an explicit reciprocity law in local class field theory via $(\varphi, \Gamma)$-modules}
\author{Naoto Dainobu}
\address{Department of Mathematics \\
3-14-1 Hiyoshi, Kohoku-ku, Yokohama-shi, Kanagawa 223-8522 Japan}
\email{vicarious@keio.jp}
\begin{document}
\maketitle
\begin{abstract}
Let $K$ be an unramified extension of $\Q_2$ and $\mu_{2^n}$ the group of $2^n$-th root of unity for a fixed integer $n \geqslant 2$. In this paper, we give an explicit formula for the $\mu_{2^n}$-valued Hilbert symbol over $K_n:=K(\mu_{2^n})$ using the theory of $(\f, \Gamma)$-modules.
\end{abstract}

\section{Introduction}
In local class field theory, we have a long tradition of describing the reciprocity map explicitly. Such a theory is usually called \textit{explicit reciprocity law}. Especially for Kummer extensions, we can study the behavior of the reciprocity map using the {\it Hilbert symbol}, which we first recall. Let $p$ be a prime number and $F$ a local field with finite residue field of characteristic $p$. Here we assume $F$ contains the group of $p^n$-th roots of unity $\mu_{p^n}$ for some $n \in\Z_{>0}$ in a fixed algebraic closure $\overline{\Q_p}$ of $\Q_p$. The Hilbert symbol over $F$ is a pairing defined as follows.

\begin{defn}[Hilbert symbol]
We define the $p^n$-th Hilbert symbol $(\cdot, \cdot)_{F, p^n}$ over $F$ as 
\[
(x, y)_{F, p^n}:=\frac{\rho_{F}(x)(\sqrt[p^n]{y})}{\sqrt[p^n]{y}}\in\mu_{p^n}\ \ (x, y \in F^{\times}),
\]
where $\rho_{F} : F^{\times} \rightarrow \Gal(F^{\mathrm{ab}}/F)$ denotes the local reciprocity map over $F$ and $F^{\mathrm{ab}}$ the maximal abelian extension of $F$.
\end{defn}

The history of explicit reciprocity law began with Kummer's work in 1858 where he essentially treated the case $F=\Q_p(\zeta_p)$ for an odd prime $p$, and gave an explicit formula for the $p$-th Hilbert symbol $(x, y)_{\Q_p(\zeta_p), p}$ for principal units $x, y$. Currently, so many types of explicit formulas are known for the Hilbert symbol. In \cite{AH}, Artin and Hasse gave such a formula of $(x, y)_{\Q_p(\zeta_{p^n}), p^n}$ for special pairs $(x, y)\in (F^{\times})^2$ as in Theorem \ref{AH} below. Iwasawa generalized their formula for more general pairs in \cite{Iwa2}, and then Coleman further generalized it in \cite{Col2}. Several generalizations of the Hilbert symbol are now also known. Wiles gave an explicit formula of the generalized Hilbert symbol for Lubin-Tate extensions of local fields in \cite{W} and de Shalit gave its generalization in \cite{dS}. The Hilbert symbol can be extended to higher local fields. Kurihara \cite{Ku} and Zinoviev \cite{Z} gave generalizations of classical Iwasawa's formula to ones for higher local fields. Fl\'orez further generalized them for an arbitrary Lubin-Tate extension in \cite{Fl}. Kato treated certain cohomological symbol defined for general local ring which is a vast generalization of the Hilbert symbol and gave an explicit formula for it in \cite{Ka}. 

Thus, the Hilbert symbol has been studied deeply by many people. However, when $p=2$, we still have a less understanding of the symbol than the case $p>2$. In fact, some formulas to compute the symbol we noted above do not work when $p=2$. For example, Kummer, Iwasawa, Wiles, de Shalit, Zinoviev, Fl\'orez and Kato's result do not work in such a case. It is because we can not apply some theory to calculate the symbol in that case. For instance, the theory of syntomic cohomology Kato used in \cite{Ka} does not work when $p=2$. Thus we often have some difficulties in the theory of explicit reciprocity law in the case $p=2$, and that is the case we treat in this paper.

In \cite{Ben}, Benois calculated the Hilbert symbol with the theory of $(\f, \Gamma)$-modules when $p$ is odd, and reproved Coleman's explicit formula. In this paper, extending this Benois' work, we give an explicit formula for the Hilbert symbol via $(\f, \Gamma)$-modules when $p=2$.

Here we describe some details of our main result. We often omit the suffix $p^n$ in the Hilbert symbol $(\cdot, \cdot)_{F, p^n}$ and write it as $(\cdot, \cdot)_{F}$ if no confusion occurs. Let $K$ be an unramified extension of $\Q_p$, $\Ok$ its ring of integers and $K_n:=K(\mu_{p^n})$. Choosing a primitive $p^n$-th root of unity $\z \in \mu_{p^n}$, we define another symbol $[\cdot, \cdot]_{K_n} : K_n^{\times} \times K_n^{\times} \rightarrow \Z/p^n$ by $(x, y)_{K_n}=\zeta_{p^n}^{[x, y]_{K_n}}$. The main result in this paper is the following formula.
\begin{thm*}[Main result]\label{main}
Suppose $n\geqslant 2$ and $p=2$. Let $U_{K_n}^1$ be the principal unit group of $K_n$. For $x, y \in U_{K_n}^1$, we have 
\begin{eqnarray*}
& &[x, y]_{K_n}\\ 
&\hspace{-3mm}=&\hspace{-2mm}-(1+2^{n-1})\Tr_{K/\Q_2}\left(\Res_{\pi_n}\left(D\log f\Log(g)-\Log(f)\f(D\log(g)\right)
\frac{d\pi_n}{\pi(1+\pi_n)}\right)\\
& &-2^n\Tr_{K/\Q_2}\left(\Res_{\pi_n}\left(\Log(f)\f(Y_y)-Y_x\Log(g)\right)\frac{d\pi_n}{\pi(1+\pi_n)}\right).
\end{eqnarray*}
Here $\pi_n$ is an indeterminate defined in Section 2, $f=f(\pi_n), g=g(\pi_n)$ are power series of $\pi_n$ in $1+\pi_n\Ok{[[\pi_n]]}$ which satisfy $f(\zeta_{p^n}-1)=x, g(\zeta_{p^n}-1)=y$, and $\Res_{\pi_n}$ denotes the residue of power series with respect to $\pi_n$. Power series $Y_x(\pi_n), Y_y(\pi_n) \in \frac{1}{2}\Ok{[[\pi_n]]}$ and operators $D, \Log$ are defined in Proposition \ref{key}. 
\end{thm*}
The first term in our formula is similar to Benois' result in \cite[Proposition 2.3.1.]{Ben}, but our formula has an extra term. It is interesting for the author to see the appearance of such an extra term since he expected that the result would be a similar one to Benois' result. We explain from where this extra term comes, describing some difficulties to extend Benois' work to the case $p=2$ and how we overcome them. 

To calculate the Hilbert symbol, Benois interpreted the Kummer map $\kappa : K_n^{\times}\rightarrow H^1(K_n, \Z_p(1))$ in terms of $(\f,\G)$-modules in \cite[Proposition 2.1.5.]{Ben}. We have an isomorphism $h^1 : H^1(K_n, \Z_p(1))\xrightarrow{\sim} H^1_{\Phi\G}(\Ak(1))$ where $H^1_{\Phi\G}(\Ak(1))$ denotes certain cohomology group defined by $(\f, \G)$-modules (see Theorem \ref{cohcomp}). For $x \in U_{K_n}^1$, Benois determined a representative of the cohomology class $h^1\circ\kappa(x)$ explicitly. This is the most essential part in his work. However, this Benois' calculation of $h^1\circ\kappa$ has $2$ in its denominator. Hence this result is no longer valid when $p=2$ since we treat cohomology groups with integral coefficients. Thus we need to calculate $h^1\circ\kappa$ with a different manner. This is the main difficulty in our case $p=2$. 

One of the main ideas to overcome this difficulty is to compute $h^1\circ\kappa$ permitting the denominators once. In other words, we use the following commutative diagram
\vspace{-2mm}
\[
\xymatrix{
 U_{K_n}^1 \ar[r]^-{\kappa} & H^1(K_n, \Z_2(1)) \ar[d]^{\iota}\ar[r]^-{h^1}_{\sim}  & H^1_{\Phi\G}(\Ak(1)) \ar[d]^{\iota_{\Phi\G}} \\
  &H^1(K_n, \Q_2(1)) \ar[r]^-{h^1_{\Q_2}}_-{\sim} &  H^1_{\Phi\G}(\Ak(1)\tens{\Z_2}{\Q_2}), 
  }\vspace{-3mm}
\]
and compute the composite homomorphism $h^1_{\Q_2}\circ\iota\circ\kappa(x)$ for $x \in U_{K_n}^1$. Here, the isomorphism $h^1_{\Q_2} : H^1(K_n, \Q_2(1)) \rightarrow H^1_{\Phi\G}(A_{K_n}(1)\otimes \Q_2)$ is a scalar extension of $h^1$ to the field of fractions. The vertical arrows $\iota$, $\iota_{\Phi\G}$ which are almost injective denote the morphisms induced by the inclusions between coefficients. We get an explicit representative of the cohomology class $h^1_{\Q_2}\circ\iota\circ\kappa(x)$ with denominators here. We do this calculation in Lemma \ref{tensorQ}, and this is the most technical part in this paper. Next, we determine a suitable new representative of the cohomology class $h^1_{\Q_2}\circ\iota\circ\kappa(x)$ explicitly within integral coefficients in the proof of Proposition \ref{key}. Then the new representative gives a cohomology class in $H^1_{\Phi\G}(\Ak(1))$, the cohomology group with integral coefficients. The image of this cohomology class under $\iota_{\Phi\G}$ is $h^1_{\Q_2}\circ\iota\circ\kappa(x)$. Thus, this new representative is exactly the one which represents $h^1\circ\kappa(x)$ due to the commutativity of the above diagram and almost injectivity of $\iota_{\Phi\G}$ (See Proposition \ref{key} for more details).

To determine a new integral representative of the cohomology class $h^1_{\Q_2}\circ\iota\circ\kappa(x)$, we subtract a suitable 1-coboundary from the old representative of $h^1_{\Q_2}\circ\iota\circ\kappa(x)$ with denominators, and make it integral. We construct such a suitable $1$-coboundary for each $x\in U_{K_n}^1$ in Lemma \ref{slide}, solving certain equation of power series. Then we show the result of the subtraction has no denominators in Lemma \ref{nodenom} using the cocycle condition of $H^1_{\Phi\G}(\Ak(1)\otimes\Q_2)$ and explicit calculations of power series.

The extra term in our formula in the main result comes from the modification of the representative of $h^1_{\Q_2}\circ\iota\circ\kappa(x)$ by subtracting the suitable $1$-coboundary. We note that our argument can yield Benois' result when $p>2$. In this case, we need no modifications of the representative of $h^1_{\Q_p} \circ \iota\circ\kappa (x)$ since $2$ is invertible, and we have no extra terms as a result.

Note also that Benois showed his result is the same as Coleman's formula in \cite{Col2}. However, because of the extra term in our formula, we do not understand precise relations between our formula and Coleman's formula for $p=2$.

From a viewpoint of the theory of $(\f, \Gamma)$-modules, the author thinks Proposition \ref{key} which is a calculation of $h^1 \circ \kappa$ is important. This is the first result which gives an interpretation of Kummer map with integral coefficients in terms of $(\f, \Gamma)$-modules when $p=2$. The author hopes Proposition \ref{key} would have some contribution to the integral theory of $(\f, \Gamma)$-modules and its application of the theory of general explicit reciprocity law of integral $p$-adic representations.

At the end of this section, we write the outline of this paper. In section 2, we introduce some basic tools such as $(\f, \Gamma)$-modules and describe how to use them for calculating the Hilbert symbol. In section 3, we give an explicit interpretation of the Kummer map in terms of $(\f, \Gamma)$-modules. Using this interpretation, we finally calculate the Hilbert symbol and show the main theorem in section 4. 
\subsection*{Acknowledgement}
The author would like to thank his supervisor Professor Masato Kurihara heartily for his continued support and helpful discussions. Thanks are also due to Professor Ivan Fesenko who gave an intensive course on class field theory in Kyoto in 2018. The course greatly led the author to this topic. The author also grateful to Professor Victor Abrashkin and Denis Benois. They kindly replied to author's questions on their paper \cite{Abs}, \cite{Ben} respectively. This research was supported by JSPS KAKENHI Grant Number 21J13502.

\section{Preliminaries}
This section is devoted to describe some fundamental tools we mainly use to compute the Hilbert symbol.
\subsection{$(\f, \G)$-modules}
 We first recall Fontaine's theory of $(\f, \G)$-modules. 
\begin{defn}
Let $\C_p$ be the $p$-adic completion of $\overline{\Q_p}$ and $\mathcal{O}_{\C_p}$ its ring of integers. We define
\[
\Ep:=\ilim \mathcal{O}_{\C_p},\ \ \E:=\ilim\C_p.
\]
Here the transition maps of projective limits are the $p$-th power homomorphisms.
\end{defn}
It is a well-known fact that $\Ep$  and $\E$ are perfect rings of characteristic $p$ under some addition defined properly and componentwise multiplication. We define a valuation $v_{\E}$ on $\E$ as $v_{\E}((x_0, x_1, \ldots)):=v_{p}(x_0)$ where $v_p$ is the $p$-adic valuation on $\C_p$ normalized as $v_p(p)=1$. Then $\Ep$ is the valuation ring of $v_{\E}$ and $\E$ is a complete discrete valuation ring with respect to $v_{\E}$. Fixing a compatible system of roots of unity $\{\z\}_{n}$ such that $\zeta^p_{p^{n+1}}=\zeta_{p^n}\ (n\geqslant 0)$, we set $\e:=(1, \zeta_p, \zeta_{p^2}, \cdots) \in \Ep$. In the following, we write $W(R)$ as the Witt ring of $R$ for a perfect ring R of characteristic $p$.
\begin{defn}
We define
\[
\Ap:=W(\Ep), \ \Ac:=W(\E).
\]
\end{defn}

Putting $\pi:=[\e]-1$, we consider the $(p, \pi)$-adic topology on $\Ap$ and $\Ac$. There is an injective map $\overline{\F} \rightarrow \Ep\ (a \mapsto ([a], [a]^{\frac{1}{p}}, [a]^{\frac{1}{p^2}}, \cdots))$ where $[\cdot]$ denotes the Teichm\"uller representative and we can identify $\overline{\F}$ as a subring of $\Ep$. Hence we can identify $\Ok$ as a subring of $\Ap$. For every integer $n\geqslant 1$, we set $\pi_n:=[\e^{\frac{1}{p^n}}]-1$ and introduce the following ring $A_{K_n}$ of power series in $\Ac$.
\begin{defn}
\[
\Ak:=\Ok\{\{\pi_n\}\}:=\left\{\sum_{m \in \mathbb{Z}} a_m \pi^m_n \mid a_m \in \Ok,\ a_m \xrightarrow[m \rightarrow -\infty]{} 0\right\}.
\]
\end{defn} 
This ring $\Ak$ is the $p$-adic completion of $\mathcal{O}_{K}((\pi_n))$. Since $\Ok((\pi_n))\subset \Ac$ and $\Ac$ is $p$-adically complete, $\Ak$ is a subring of $\Ac$. We put $A_n$ as the $p$-adic completion of the maximal unramified extension of $\Ak$ in $\Ac$. Let $K_{\cyc}:=K(\zeta_{p^{\infty}})$ and $\Gn:=\Gal(K_{\cyc}/K_n)$. We assume $\Gn$ is a procyclic group. When $p=2$, this holds if $n\geqslant 2$ while this holds automatically when $p$ is odd. We fix a topological generator $\g$ of $\Gn$. Here we see actions of $\Gn$ and Frobenius $\f$ on $A_{K_n}$. Since there is a componentwise action of $G_{K_n}$ on $\E$, we have an action of $G_{K_n}$ on its Witt ring $\Ac$. This action is stable on the subring $A_n$ and it is well-known that $A_n^{G_{K_{\cyc}}}=\Ak$. Thus the quotient group $\Gn=G_{K_n}/G_{K_{\cyc}}$ acts on $\Ak$. We can see that $\g$ acts on $\pi_n$ as $\g(\pi_n)=(1+\pi_n)^{\ch(\g)}-1$ and on the coefficient ring $\Ok$ trivially, where $\ch$ denotes the $p$-adic cyclotomic character. On the other hand, we have the Frobenius homomorphism $\f$ on $\Ac=W(\E)$ as the lift of $p$-th power homomorphism on $\E$. This induces an action of $\f$ on the subring $\Ak \subset \Ac$. We can see that $\f$ acts on $\pi_n$ as $\f(\pi_n)=(1+\pi_n)^p-1$ and on the coefficient ring $\Ok$ as the Frobenius element in $\Gal(K/\Q_p)$.
\begin{defn}\label{fg-mod}
A $(\f, \Gn)$-module over $\Ak$ is a finitely generated $\Ak$-module equipped with continuous semilinear actions of $\f$ and $\Gn$ which commute with each other.
\end{defn}
Let $\Bp:=\Ap\tens{\Z_p}\Q_p$, $\Bc:=\Ac\tens{\Z_p}\Q_p$, $\Bk:=\Ak \tens{\Z_p}\Q_p$ and $B_n:=A_n \tens{\Z_p}\Q_p$. We can also define the notion of $(\f, \Gn)$-modules over $\Bk$ in the same way as Definition \ref{fg-mod}.

\subsection{$p$-adic representations and $(\f, \G)$-modules}
In \cite{F}, Fontaine proved the following striking theorem.
\begin{thm}[Fontaine]\label{cat}
Let $\mathrm{Rep}_{\Z_p} G_{K_n}$ be the category of $p$-adic representations of $G_{K_n}$ over $\Z_p$ and $\Phi\Gamma^{\text{\'et}}_{\Ak}$ the category of \'etale $(\f, \Gn)$-modules over $\Ak$. Then there is a category equivalence
\[
\mathbf{D}: \mathrm{Rep}_{\Z_p}\Gkn \overset{\sim}{\longrightarrow} \Phi\Gamma^{\text{\'et}}_{\Ak},
\]
where for an object $T$ in $\mathrm{Rep}_{\Z_p} \Gkn$, the functor $\mathbf{D}$ is defined as 
\[
\mathbf{D}(T)=(T\underset{\Z_p}{\otimes}A_n)^{G_{K_{\cyc}}}.
\]
Here, we consider a diagonal action of $\Gn$ and an action of $\f$ only on the right component $A_n$ on $\mathbf{D}(T)$.
\end{thm} 
We do not define the notion of \'etale $(\f, \Gn)$-module. Here is an example of Theorem \ref{cat}. Let $T:=\Z_p(1):=\ilim \mu_{p^n}$, then 
\[
\mathbf{D}(\Z_p(1))=(\Z_p(1)\tens{\Z_p}A_n)^{G_{K_{\cyc}}}=(A_n(1))^{G_{K_{\cyc}}}=\Ak(1).
\]
In the above computation, we define $\Ak(1):=\Z_p(1)\tens{\Z_p}\Ak$.

The similar category equivalence exists between the category $\mathrm{Rep}_{\Q_p} \Gkn$ of $p$-adic representations over $\Q_p$ and the category $\Phi\Gamma^{\textit{\'et}}_{\Bk}$ of \'etale $(\f, \Gn)$-modules over $\Bk$.
\begin{thm}[Fontaine]\label{catp}
There is a category equivalence
\[
\mathbf{D}: \mathrm{Rep}_{\Q_p} \Gkn \overset{\sim}{\longrightarrow} \Phi\Gamma^{et}_{\Bk}
\]
where $\mathbf{D}(V):=(V\underset{\Q_p}{\otimes}B_n)^{G_{K_{\cyc}}}$ for an object $V$ in $\mathrm{Rep}_{\Z_p} \Gkn$.
\end{thm}
We can compute the Galois cohomology group of $T \in \mathrm{Rep}_{\Z_p} \Gkn$ using the corresponding $(\f, \Gn)$-module $\mathbf{D}(T)$.
\begin{defn}[Fontaine-Herr]\label{cpx}
Let $T$ be an object in $\mathrm{Rep}_{\Z_p} \Gkn$. For the corresponding $(\f, \Gn)$-module $\mathbf{D}(T)$, we define a complex
\[
C^{\bullet}(\mathbf{D}(T)): 0\rightarrow \mathbf{D}(T) \xrightarrow[\alpha]{} \mathbf{D}(T)^{\oplus 2} \xrightarrow[\beta]{} \mathbf{D}(T) \rightarrow 0,
\]
where the maps $\alpha, \beta$ defined as
\begin{eqnarray*}
\alpha(z)&:=&[\left((\f-1)(x), (\g-1)(x)\right)]\ \ (x\in \mathbf{D}(T)),\\
\beta(y, z)&:=& [(\g-1)(y)+(1-\f)(z)]\ \ (y, z \in \mathbf{D}(T)).
\end{eqnarray*}
\end{defn}
In the following, we write the cohomology group $H^i(C^{\bullet}(\mathbf{D}(T)))$ as $H^i_{\Phi\Gamma}(\mathbf{D}(T))$.
\begin{thm}[Fontaine-Herr]\label{cohcomp}
Let $T$ be an object in $\mathrm{Rep}_{\Z_p} \Gkn$. For each $i\geqslant 0$, we have an isomorphism
\[
h^i : H^i(K_n, T) \xrightarrow{\sim} H^i_{\Phi\Gamma}(\mathbf{D}(T)).
\]
\end{thm}

Thanks to Theorem \ref{cohcomp}, for example, an element in $H^1(K_n, \Z_p(1))$ correspond to a cohomology class in $H^1_{\Phi\Gamma}(\Ak(1))$ represented by a pair of power series in $\Ak(1)$ via $h^1$. In the succeeding sections, we use this explicit interpretation of Galois cohomology classes to compute the Hilbert symbol. 

We note that exactly the same statement as Theorem \ref{cohcomp} holds for $p$-adic representation $V$ over $\mathbb{Q}_p$. 
\begin{thm}[Fontaine-Herr]\label{cohcompp}
Let $V$ be an object in $\mathrm{Rep}_{\Q_p} \Gkn$. For each $i\geqslant 0$, we have an isomorphism
\[
h^i_{\mathbb{Q}_p} : H^i(K_n, V) \xrightarrow{\sim} H^i_{\Phi\Gamma}(\mathbf{D}(V)):=H^i(C^{\bullet}(\mathbf{D}(V))).
\]
Here, the complex $C^{\bullet}(\mathbf{D}(V))$ of $(\f, \Gn)$-modules over $\Bk$ defined the same way as in Definition \ref{cpx}.
\end{thm}
We can compute a cup product of Galois cohomology groups using that of $(\f, \Gn)$-modules and isomorphism $h^i$.
\begin{prop}[Fontaine-Herr]\label{cup}
Let $T_1, T_2$ be objects in $\mathrm{Rep}_{\Q_p} \Gkn$. We define a bilinear pairing $\cup_{\Phi\G} : H^1_{\Phi\G}(\mathbf{D}(T_1))\times H^1_{\Phi\G}(\mathbf{D}(T_2)) \rightarrow H^2_{\Phi\G}(\mathbf{D}(T_1 \otimes T_2))$ as 
\[
[(m_1, n_1)] \cup_{\Phi\G} [(m_2, n_2)] := [n_1 \otimes \g(m_2)-m_1 \otimes\f(n_2)],
\]
where $m_1, n_1 \in \mathbf{D}(T_1)$ and $m_2, n_2 \in \mathbf{D}(T_2)$. Then the following diagram is commutative.
\[
  \xymatrix{
  H^1(K_n, T_1) \times H^1(K_n, T_2) \ar[r]^{\hspace{6mm}\cup}\ar[d]_{h^1 \times h^1} & H^2(K_n, T_1 \otimes T_2)\ar[d]^{h^2} \\
  H^1_{\Phi\G}(\mathbf{D}(T_1)) \times H^1_{\Phi\G}(\mathbf{D}(T_2)) \ar[r]^{\hspace{10mm}\cup_{\Phi\G}} &  H^2_{\Phi\G}(\mathbf{D}(T_1 \otimes T_2)) 
  }
\]
\end{prop}
Note that Fontaine and Herr gave cup products of cohomology groups of $(\f, \Gn)$-modules for other degrees than $H^1$. See \cite{H1} or \cite{H2} for details.

Finally, we introduce an isomorphism $\mathrm{TR}_{K_n} : H^2_{\Phi\G}(\Ak(1)) \rightarrow \Z_p$ corresponding the invariant map $\mathrm{inv}_{K_n} : H^1(K_n, \Z_p(1)) \rightarrow \Z_p$ in local class field theory. In the following, we consider $\e$ as a basis of the Tate twist $\Z_p(1)$ and write $a \otimes \e$ for $a \in \Ak$ when we consider $a$ as an element in $\Ak(1)$. In \cite{Ben}, Benois proved the following result.
\begin{prop}[Benois]\label{trace}
Define $\mathrm{TR}_{K_n} : H^2_{\Phi\G}(\Ak(1)) \rightarrow \Z_p$ as 
\[
\mathrm{TR}_{K_n}([a\otimes\e]):= -\frac{p^n}{\log(\ch(\g))}\mathrm{Tr}_{K/\Q_p}\Res_{\pi_n} \left(\frac{a d\pi_n}{1+\pi_n}\right) \ \ (a \in \Ak),
\]
where for an element $f(\pi_n)d\pi_n = \left(\sum_{i \in \Z} a_i \pi_n^i \right)d\pi_n$ of an $\Ok$-module of differential 1-forms $\Omega_{\Ak/\Ok}^1$, we define $\Res(f(\pi_n)):=a_{-1}$. Then the following diagram is commutative:
\[
  \xymatrix{
    H^2(K_n, \Z_p(1))\ar[d]^{h^2}\ar[r]^{\hspace{10mm}\mathrm{inv}_{K_n}}&\Z_p\\
     H^2_{\Phi\G}(\Ak(1))\ar[ru]_{\mathrm{TR}_{K_n}}& 
  }
\]
\end{prop}
\begin{rem}
Although Benois proved the above result for an odd prime $p$, we can check the result is also valid for $p=2$ by the similar way in \cite{Ben}.
\end{rem}
\subsection{Fontaine's crystalline period ring}\label{crys} In our calculation of the Hilbert symbol, we use Fontaine's crystalline period ring $\A$ which we recall below.
\begin{defn}
We define a ring homomorphism $\theta$ as 
\[
\theta: \Ap \rightarrow \mathcal{O}_{\C_p},\ \ \ \ \sum^{\infty}_{i=0}[x_i] p^i \mapsto \sum^{\infty}_{i=0}(x_i)_0 p^i 
\]
where $x_i\in \Ep$ and $(x_i)_0\in \mathcal{O}_{\C_p}$ denotes its 0-th component.
\end{defn}
This is a homomorphism of $\mathcal{O}_{\Q^{\mathrm{ur}}_p}$-algebra where $\Q^{\mathrm{ur}}_p$ denotes the maximal unramified extension of $\Q_p$ and $\mathcal{O}_{\Q^{\mathrm{ur}}_p}$ its ring of integers. Put $v:=\pi/\pi_1=1+[\e^{\frac{1}{p}}]+[\e^{\frac{1}{p}}]^2+\cdots+[\e^{\frac{1}{p}}]^{p-1}$. Then it is a well-known fact that the kernel of $\theta$ is principal and generated by $v$. We put $A^0_{\mathrm{crys}}:=\Ap{[\{\frac{v^m}{m!}\}_{m>0}]}$, the divided power envelop of $\Ap$ with respect to $\mathrm{Ker} \theta$. We define $\A$ as its $p$-adic completion. More explicitly,  
\[
\A=\left\{\sum^{\infty}_{m=0}a_m \frac{v^m}{m!}\ \ \middle|\ \ a_m \rightarrow 0\ (m \rightarrow \infty)\ \text{$p$-adically}\right\}.
\]
We define an element $t \in \A$ as
\[
t:=\log(1+\pi)=\sum_{m=1}^{\infty} (-1)^{m+1}\frac{\pi^m}{m}.
\]
In fact, this infinite sum converges in $\A$ with respect to its $p$-adic topology. We put $\B:=A_{\mathrm{crys}}[\frac{1}{p}]$ and $B_{\mathrm{crys}}:=\B{[\frac{1}{t}]}$. Here we state a lemma we use in the next section.
\begin{lem}\label{logtoreru}
Suppose $a \in \Ap$ satisfies $\theta(a)=1$, then 
\[
\log a := \sum_{m=1}^{\infty} (-1)^{m+1}\frac{(a-1)^m}{m}
\] 
converges in $\A$.
\end{lem}
(Proof of Lemma \ref{logtoreru})

Since $\theta(a)=1$, there exist $x \in \Ap$ such that $a=1+xv$. Then we have,
\[
\log a = \log (1+xv) = \sum_{m=1}^{\infty} (-1)^{m+1}\frac{(xv)^m}{m}.
\]
While
\[
(-1)^{m+1}\frac{(xv)^m}{m}=(-1)^{m+1}(m-1)!\cdot x^m\cdot \frac{v^m}{m!}.
\]
The factor $(m-1)!$ converges to $0$ as $m \rightarrow \infty$ with respect to the $p$-adic topology in $\A$, which implies the convergence of $\log a$ in $\A$.\hfill$\square$
\subsection{Strategy of the calculation}\label{str}In this subsection, we briefly describe the method of calculation. We mainly follow Benois' strategy in \cite{Ben}.　There is an exact sequence of $\Gkn$-modules
\[
1 \rightarrow \mu_{p^m} \rightarrow \overline{K_n} \rightarrow \overline{K_n} \rightarrow 1
\]
from which we get $\kappa_m : K^{\times}_n \rightarrow H^1(K_n, \mu_{p^m})
$ as its connecting homomorphism. Taking the inverse limit with respect to $m$, we have 
\[
\ka : K^{\times}_n \rightarrow H^1(K_n, \Z_p(1))
\]
which we call the Kummer map. Using this $\kappa$, we have the following cohomological interpretation of the Hilbert symbol.
\[
  \xymatrix{
  (K_n^{\times})^{\otimes 2} \ar[d]^{\kappa^{\otimes 2}}\ar@/^20pt/[drrr]_{(\cdot, \cdot)_{K_n}}  & & & \\
  H^1(K_n, \Z_p(1))^{\otimes 2} \ar[r]^{\cup_{\Gal}}  & H^2(K_n, \Z_p(2))\ar[r]^-{\mod p^n} & H^2(K_n, \mu_{p^n})\otimes\mu_{p^n}\ar[r]^-{\hspace{10mm}\mathrm{inv}_{K_n}}_-{\sim}&\mu_{p^n},
    }
\]
where $\cup_{\Gal}$ denotes the cup product of Galois cohomology groups. Note that since $K_n$ contains $\mu_{p^n}$, we have an isomorphism $H^2(K_n, \mu_{p^n}^{\otimes 2}) \xrightarrow{\sim} H^2(K_n, \mu_{p^n})\otimes \mu_{p^n}$ which induced by the cup product. On the other hand, the morphisms in the second row can be calculated using the theory of $(\f, \G)$-modules as
\[
  \xymatrix{
  H^1(K_n, \Z_p(1))^{\otimes 2} \ar[d]^{{h^1}^{\otimes 2}}\ar[r]^{\cup_{\Gal}}  & H^2(K_n, \Z_p(2)) \ar[d]^{h^2}\ar[r]^{\hspace{-5mm}\mod p^n} & H^2(K_n, \mu_{p^n})\otimes\mu_{p^n}\ar[d]^{h^2}\ar[r]^{\hspace{13mm}\mathrm{inv}}_{\hspace{13mm}\sim}&\mu_{p^n}\\
  H^1_{\Phi\G}(\Ak(1))^{\otimes 2} \ar[r]^{\cup_{\Phi\G}} &  H^2_{\Phi\G}(\Ak(2))\ar[r]^{\hspace{-5mm}\mod p^n} & H^2_{\Phi\G}(\mu_{p^n})\otimes\mu_{p^n}\ar[ru]_{\overline{\mathrm{TR}_{K_n}}},& 
  }
\]
where $\cup_{\Phi\G}$ is the cup product we define in Proposition \ref{cup} and $\overline{\mathrm{TR}_{K_n}}$ is the mod $p^n$ reduction of the isomorphism $\mathrm{TR}_{K_n}$ in Proposition \ref{trace}. Since $\cup_{\Phi\G}$ and $\overline{\mathrm{TR}_{K_n}}$ are given explicitly, all we have to do for the calculation of the Hilbert symbol is an explicit computation of the composite homomorphism $h^1 \circ \kappa$. 
\begin{rem}
Kato computed this cup product $\cup_{\Gal}$ via the theory of syntomic cohomology in \cite{Ka} for more general setting when the residue characteristic $p$ is odd. Note that this cohomology theory does not work for our case $p=2$.
\end{rem} 
\section{Calculation of the Kummer map}
In this section, we compute the composite homomorphism $h^1 \circ \kappa$.
\subsection{explicit calculation of the isomorphism $h^1$} First, we give an explicit formula of the isomorphisms 
\[h^1 : H^1(K_n, \Z_p(1)) \rightarrow H^1_{\Phi\G}(\Ak(1)),\ \ h^1_{\Q_p(1)} : H^1(K_n, \Q_p(1)) \rightarrow H^1_{\Phi\G}(\Bk(1)).
\]
\begin{prop}\label{isocomp}
For a cohomology class $[c]\in H^1(K_n, \Z_p(1))$ $(resp.\ H^1(K_n, \Q_p(1)))$ which is represented by a 1-cocycle $c : \Gkn \rightarrow \Z_p(1)\ (resp.\ \Q_p(1)) ,\ g \mapsto c(g)\otimes\e$, we have 
\[
h^1([c])=\left[(\f-1)(\xi_{c}\otimes \e), (\widehat{\g}-1)(\xi_{c}\otimes \e)+c(\widehat{\g})\otimes\e  \right].
\]
\[
(resp.\ h^1_{\Q_p(1)}([c])=\left[(\f-1)(\xi_{c}\otimes \e), (\widehat{\g}-1)(\xi_{c}\otimes \e)+c(\widehat{\g})\otimes\e\right].)
\]
Here, $\widehat{\g}$ is any lift of $\g$ to $\Gkn$ and $\xi_{c} \in A_n$ $(resp.\ B_n)$ is an element which satisfies
\[
g(\xi_{c})=\xi_{c}-c(g) \ (\forall g \in G_{K_{\cyc}}).
\]
\end{prop}
(Proof of Proposition \ref{isocomp})

Since computations for $h^1$ and $h^1_{\Q_p(1)}$ are exactly the same, we give a proof only for $h^1$. The cohomology class $[c]\in H^1(K_n, \Z_p(1))$ corresponds to the following extension of $\Z_p$ by $\Z_p(1)$ as a $\Gkn$-module:
\[
0\rightarrow \Z_p(1) \rightarrow T_{[c]} \xrightarrow{f} \Z_p \rightarrow 0
\]
We take $1\otimes \e$ and $e$ as a basis of $T_{[c]}$ over $\Z_p$ where $g \in \Gkn$ acts on $e$ as $g(e)= e+c(g)\otimes\e$. Then for an element $x:=a \otimes \e + b\cdot e \in T_{[c]}\ (a, b \in \Z_p)$, the homomorphism $f$ is given by $f(x)=b$. Applying the functor $\mathbf{D}$, which is exact, we have a corresponding exact sequence of $(\f, \Gn)$-modules 
\[
0\rightarrow \Ak(1) \rightarrow \mathbf{D}(T_{[c]}) \xrightarrow{\mathbf{D}(f)} \Ak \rightarrow 0.
\]
Putting $\delta : \Z_p \rightarrow H^1(K_n, \Z_p(1))$ and $\delta_{\Phi\G} : \Z_p=H^0_{\Phi\G}(\Ak) \rightarrow H^1_{\Phi\G}(\Ak(1))$ as the connecting homomorphisms of the above exact sequences respectively, we have a commutative diagram
\[
  \xymatrix{
  \Z_p \ar[d]^{h^0 = id}\ar[r]^{\hspace{-8mm}\delta}  & H^1(K_n, \Z_p(1)) \ar[d]^{h^1}\\
  H^0_{\Phi\G}(\Ak)=\Z_p \ar[r]^{\delta_{\Phi\G}} &  H^1_{\Phi\G}(\Ak(1)). 
  }
\]
Since $\delta(1)=[c]$, we know $\delta_{\Phi\G}(1)=h^1([c])$ by the above diagram. So we compute $\delta_{\Phi\G}(1)$ following the definition of the connecting homomorphism. By the definition of the functor $\mathbf{D}$ we have, 
\begin{eqnarray*}
\mathbf{D}(T_{[c]})=\left(T_{[c]} \otimes A_n\right)^{G_{K_{\cyc}}}& = &\left(\Z_p(1) \oplus \Z_p\cdot e\right)^{G_{K_{\cyc}}}\\
& = & \left( A_n(1) \oplus A_n \cdot e \right)^{G_{K_{\cyc}}}.
\end{eqnarray*}
For an element $x:=a\otimes \e + b\cdot e\in A_n(1) \oplus A_n \cdot e\ (a, b \in A_n)$ and $g \in G_{K_{\cyc}}$, 
\begin{eqnarray*}
g(x)=g(a\otimes \e + b\cdot e)& = & \ch(g)g(a)\otimes \e + g(b)(e+c(g)\otimes\e)\\
&=& (g(a)+g(b)c(g))\otimes\e+g(b)\cdot e.
\end{eqnarray*}
Thus $x=a\otimes \e + b\cdot e$ is fixed by $G_{K_{\cyc}}$ if and only if 
\[
g(a)+g(b)c(g)=a,\ g(b)=b \ \ (\forall g \in G_{K_{\cyc}}).
\]
From the second condition, $b \in (A_n)^{G_{K_{\cyc}}}=\Ak$ and thus the first condition says $g(a)+bc(g)=a\ (\forall g \in G_{K_{\cyc}})$. Hence 
\[
\mathbf{D}(T_{[c]})=\left\{a\otimes \e + b\cdot e \mid a \in A_n, b \in \Ak\ s,t\ \ g(a)+bc(g)=a\ (\forall g \in G_{K_{\cyc}})\right\}.
\]
Now we compute $\delta_{\Phi\G}(1)$. First we pick $\xi_c\otimes\e+e\in\mathbf{D}(T_{[c]})$ for some $\xi_c \in A_n$ satisfying $g(\xi_c)=\xi_c-c(g)$ for all $g \in G_{K_{\cyc}}$. This element maps to $1 \in \Ak$ under $\mathbf{D}(f)$ and we compute its image under the homomorphism $\alpha$ in Definition \ref{cpx} as
\[
\alpha(\xi_c\otimes\e+e) = ((\f-1)(\xi_c \otimes \e + e), (\g-1)(\xi_c \otimes\e + e)).
\]
On the first component, we have
\[
(\f-1)(\xi_c \otimes \e + e)=(\f(\xi_c) \otimes \e +e)-(\xi_c\otimes \e +e)=(\f-1)(\xi_x)\otimes\e,
\]
and the second component,
\begin{eqnarray*}
(\g-1)(\xi_c \otimes\e + e) &=& (\gh-1)(\xi_c\otimes \e)+(\gh-1)(e) \\
&=& (\gh-1)(\xi_c\otimes \e) + c(\gh)\otimes\e.
\end{eqnarray*}
Here since each term $\xi_c\otimes\e$ and $e$ respectively are not fixed by $G_{K_{\cyc}}$ although the element $\xi_c\otimes\e+e$ is fixed by $G_{K_{\cyc}}$, we have to take some extension $\gh$ of $\g\in\Gn$ to $G_{K_n}$ in the above computation. Thus we get
\[
\alpha(\xi_c\otimes\e+e)=\left((\f-1)(\xi_x)\otimes\e, (\gh-1)(\xi_c\otimes \e) + c(\gh)\otimes\e\right).
\]
The cohomology class in $H^1_{\Phi\G}(\Ak(1))$ defined by this pair is nothing other than the image $\delta_{\Phi\G}(1)$ by the definition of the connecting homomorphism. Hence we obtain the proposition.\hfill$\square$

\subsection{Computation of $h^1\circ\ka$}
In the following, we set $p=2$. This subsection is devoted to the computation of the homomorphism
\[
h^1\circ \kappa : K^{\times}_n \rightarrow H^1(K_n, \Z_2(1)) \rightarrow H^1_{\Phi\G}(\Ak(1)).
\]
We put $(U^1_{K_n})^{f}$ as the free part of the principal unit group $U^1_{K_n}=\langle \z \rangle \oplus (U^1_{K_n})^{f}$ of $K_n$ as a $\Z_2$-module. The following is a key proposition for our main result.
\begin{prop}\label{key}
For $x \in (U^1_{K_n})^{f}$, we have 
\[
h^1\circ\ka (x)=\left[\Log(f(\pi_n))\cdot\frac{1}{\pi} \otimes \e,\ \lambda_x(\pi_n) \otimes\e -(\ch(\g)-1)Y_x(\pi_n)\otimes \e\right],
\] 
where $f(\pi) \in 1+\pi_n\Ok{[[\pi_n]]}$ is a power series which satisfies $f(\zeta_{2^n}-1)=x$ for which the operator $\Log$ defined as $\Log(f(\pi_n)):=(\frac{\f}{p}-1)\log(f(\pi_n))$. The power series $\lambda(\pi_n) \in \Ok{[[\pi_n]]}$ is uniquely determined one corresponding to the 1-st component and satisfies
\[
\lambda_x(\pi_n) \equiv \frac{\ch(\g)-1}{2^n}D\log f(\pi_n)\ \ \mathrm{mod}\ \pi\Ok{[[\pi_n]]},
\]
where $D:=(1+\pi_n)\frac{d}{d\pi_n}$. The power series $Y_x(\pi_n) \in \frac{1}{2}\Ok{[[\pi_n]]}$ is defined as 
\[
Y_x(\pi_n):=\frac{1}{2}\sum_{i=0}^{\infty} \f^i(\Log(f(\pi_n))).
\]
\end{prop}
Although the power series $Y_x(\pi_n)$ itself has a denominator, the term $(\ch(\g)-1)Y_x(\pi_n)$ in the second component of $h^1\circ \ka(x)$ is an element of $\Ak$ since $\ch(\g)-1 \in 2^n\Z_2\ (n\geqslant 2)$.  We prove this key proposition after introducing some lemmas. First we consider a situation tensored with $\Q_2$, in other words, we think $\ka(x) \in H^1(K_n, \Z_p(1))$ as an element of $H^1(K_n, \Q_p(1))$ and compute the image of $\ka(x)$ under the isomorphism $h^1_{\Q_2}: H^1(K_n, \Q_2(1)) \rightarrow H^1_{\Phi\Gamma}(\Bk(1))$ in Theorem \ref{cohcompp}.

\begin{lem}\label{tensorQ}
For $x \in (U^1_{K_n})^{f}$,
\[
h^1_{\Q_2}\circ\kappa(x)= \left[\Log(f(\pi_n))\cdot\left(\frac{1}{2}+\frac{1}{\pi} \right)\otimes \e,\ \lambda_x (\pi_n) \otimes \e \right],
\]
where $f(\pi_n)$ is the same as Proposition \ref{key} and $\lambda_x (\pi_n)\in \Ok{[[\pi_n]]}\tens{\Z_2}\Q_2$ satisfies
\[
\lambda_x(\pi_n) \equiv \frac{\chi(\g)-1}{2^n}D\log f(\pi_n)\ \ \mathrm{mod}\ \pi\Ok{[[\pi_n]]}\tens{\mathbb{Z}_2}\mathbb{Q}_2.
\]
\end{lem}
(Proof of Lemma \ref{tensorQ})

 From Proposition \ref{isocomp}, it suffices to construct $\xi_{\ka(x)} \in B_n$ explicitly and compute actions of $\f$ and $\widehat{\g}$ on it. Put $\omega_x:=[x, x^{\frac{1}{p}}, x^{\frac{1}{p^2}}, \ldots] \in \Ep$ and  $a_x:=\frac{f(\pi_n)}{[\omega_x]} \in \Ap$. Applying $\theta: \Ap \rightarrow \mathcal{O}_{\C_p}$ which we defined in subsection \ref{crys} on $a_x$, we have
  \[
 \theta(a_x)=\theta\left(\frac{f(\pi_n)}{[\omega_x]}\right)=\frac{f(\theta(\pi_n))}{x}=\frac{f(\zeta_{2^n}-1)}{x}=1.
 \]
 Thus $\log a_x$ defines a well-defined element in $\A$ from Lemma \ref{logtoreru}. Since $\theta(a_x)=1$, there exists $a \in \Ap$ such that $a_x=1+a v$ and $\log (a_x)$ can be expressed as
\begin{eqnarray}\label{loga}
\log a_x=a v - \frac{(a v)^2}{2} + \frac{(a v)^3}{3} - \cdots +(-1)^{m+1}\frac{(a v)^m}{m} + \cdots.
\end{eqnarray}
\begin{sublem}\label{modsuru}
There exists an element $b_x \in \Ap$ such that 
\[
b_x \equiv \log a_x+\frac{\pi}{2} a^2 \mod \pi^2\B.
\]
\end{sublem}
(Proof of Sublemma \ref{modsuru}.)

Since $\E$ has characteristic 2, we have 
\[
\left(\frac{\e-1}{\e^{1/2}-1}\right)^2=(\e^{1/2}-1)^2=\e-1.
\]
Takeing the Teichm\"uller lift of the both sides, we obtain $\left[\frac{\e-1}{\e^{1/2}-1}\right]^2=[\e-1]$, and hence
 \[
v^2=\left(\frac{[\e]-1}{[\e^{1/2}]-1}\right)^2 \equiv \left[\frac{\e-1}{\e^{1/2}-1}\right]^2=[\e-1]\equiv \pi \mod 2\Ap.
\]
Thus there exists $\alpha \in \Ap$ such that $v^2 = \pi+2\alpha$. We show that  the $m$-th term $(-1)^{m+1}\frac{(a v)^m}{m}$ in (\ref{loga}) has a suitable representative $c_m$ in $\Ap$ when considered with mod $\pi^2\B$ for every $m> 2$.
\vspace{1mm}
 
(Case 1 : $2\nmid m$) 

In this case, $(-1)^{m+1}\frac{(a v)^m}{m} \in \Ap$ and we see that
\[
v^m=\frac{\pi^m}{\pi_1^m}=\frac{(\pi_1^2+2\pi_1)^m}{\pi_1^m}=(\pi_1+2)^m.
\] 
This converges to 0 as $m \rightarrow \infty$ in $\Ap$ and so does $c_m:=(-1)^{m+1}\frac{(a v)^m}{m}$.
\vspace{1mm}

(Case 2 : $2 \mid m$ and $m > 2$)\vspace{1mm}

Writing $m=2^{\ell}\cdot s$ ($2\nmid s, \ell\geqslant1$), we have
\[
(-1)^{m+1}\frac{(a v)^m}{m} = \frac{(-1)^{m+1}}{s} a^m \cdot \frac{(v^2)^{2^{\ell-1}s}}{2^{\ell}} = \frac{(-1)^{m+1}}{s} a^m \cdot \frac{(\pi + 2\alpha)^{2^{\ell-1}s}}{2^{\ell}}.
\]
On the last factor, we see that 
\begin{eqnarray*}
\frac{(\pi + 2\alpha)^{2^{\ell-1}s}}{2^{\ell}}&=&\left(\frac{\pi}{2}+\alpha\right)^{2^{\ell-1}s}\cdot 2^{2^{\ell-1}s-\ell}\\
&\equiv& \alpha^{2^{\ell-1}s} 2^{2^{\ell-1}s-\ell} + 2^{2^{\ell-1}s-2}s \alpha^{2^{\ell-1}s-1} \pi\ \mod \pi^2\B,
\end{eqnarray*}
where since $m>2$, we have $2^{\ell-1}s-2\geqslant 0$. The right-hand side converges when $m\rightarrow \infty$. We put 
\[
c_m:=\frac{(-1)^{m+1}}{s} a^m\cdot \left(\alpha^{2^{\ell-1}s} 2^{2^{\ell-1}s-\ell} + 2^{2^{\ell-1}s-2}s\alpha^{2^{\ell-1}s-1}\pi\right).
\] 

Then $c_m \in \Ap$ is congruent to $(-1)^{m+1}\frac{(a v)^m}{m}$ mod $\pi^2\B$ and converges to 0 as $m \rightarrow \infty$.

Finally, on the second term in (\ref{loga}), we see that 
\[
(-1)^{2+1}\frac{(a v)^2}{2}=-\frac{a^2(\pi+2\alpha)}{2}=-\frac{\pi}{2}a^2-\alpha.
\]
 
Then we obtain 
\[
\log a_x\equiv av-\frac{\pi}{2}a^2-\alpha+\sum_{m\geqslant3}^{\infty}c_m \mod \pi^2\B
\] 
This implies that $\log a_x+\frac{\pi}{2}a^2 \mod \pi^2\B$ is represented by a well-defined element $b_x:=av-\alpha+\sum_{m\geqslant3}^{\infty}c_m \in \Ap$. \hfill$\square$\vspace{2mm}

We go back to the proof of Lemma \ref{tensorQ}. First, we consider $G_{K_{\cyc}}$ action on this element $b_x \in \Ap$.
\begin{sublem}\label{galactb}
For $g \in G_{K_{\cyc}}$, 
\[
g(b_x) \equiv b_x-\kappa(x)(g)\pi \mod \pi_1\pi\B
\]
\end{sublem}
(Proof of sublemme \ref{galactb}.)

For $g \in G_{K_{\cyc}}$, we have
 \begin{eqnarray*}
 g(\log a_x)=\log\frac{g(f(\pi_n))}{[g(\omega_x)]}&=&\log\frac{f(\pi_n)}{[\omega_x][\e]^{\kappa(x)(g)}}\\
 &=&\log\frac{f(\pi_n)}{[\omega_x]}-\kappa(x)(g)t\nonumber\\
 &=&\log a_x-\kappa(x)(g)t.\nonumber
 \end{eqnarray*}
Here, the element $t$ is the one we defined in subsection \ref{crys}. From Sublemma \ref{modsuru}, this implies a congruence
\begin{eqnarray}\label{gactb}
g\left(b_x-\frac{\pi}{2}a^2\right) \equiv b_x-\frac{\pi}{2}a^2-\kappa(x)(g)\pi \mod \pi_1\pi\B.
\end{eqnarray}
Note that we use a congruence of mod $\pi_1\pi\B$ here which is immediately deduced from Sublemma \ref{modsuru}. Since $a=\left(\frac{f(\pi_n)}{[\omega_x]}-1\right)\cdot\frac{1}{v}$, we have  
\begin{eqnarray*}
g(a)=\left(\frac{f(\pi_n)}{[\omega_x][\e]^{\kappa(x)(g)}}-1\right)\cdot\frac{1}{v}
&=& \left(\frac{f(\pi_n)}{[\omega_x](1+\pi)^{\kappa(x)(g)}}-1\right)\cdot\frac{1}{v}\\
&\equiv& \left(\frac{f(\pi_n)}{[\omega_x]}-1\right)\cdot\frac{1}{v}=a \mod \pi_1\B.
\end{eqnarray*}
Hence we see that
\[
g\left(\frac{\pi}{2}a^2\right)=\frac{\pi}{2}g(a)^2\equiv \frac{\pi}{2}a^2 \mod \pi_1\pi\B.
\]
From (\ref{gactb}), this congruence yields
\[
g(b_x) \equiv b_x-\kappa(x)(g)\pi \mod \pi_1\pi\B. 
\]
\hfill$\square$

Next, we consider the action of $\f$ on $b_x$.
\begin{sublem}\label{factb}
\[
\left(\frac{\f}{2}-1\right)b_x\equiv \Log(f(\pi_n))+\frac{\pi}{2}(\f-1)(a^2) \mod \pi_1\pi\B
\]
\end{sublem}

(Proof of Lemma \ref{factb})

On the action of $\f$ on $\log a_x$, we see that 
\begin{eqnarray*}
\left(\frac{\f}{2}-1\right)\log a_x&=&\left(\frac{\f}{2}-1\right)\log\frac{f(\pi_n)}{[\omega_x]}\\
&=& \frac{1}{2}\log \frac{\f(f(\pi_n))}{[\omega_x]^2}-\log\frac{f(\pi_n)}{[\omega_x]}\\
&=& \left(\frac{\f}{2}-1\right)\log f(\pi_n)=\Log(f(\pi_n)).
\end{eqnarray*}
On the other hand, 
\begin{eqnarray*}
\left(\frac{\f}{2}-1\right)\frac{\pi}{2}a^2=\frac{1}{4}\f(\pi)\f(a^2)-\frac{\pi}{2}a^2&=&\frac{1}{4}(\pi^2+2\pi)\f(a^2)-\frac{\pi}{2}\\
&\equiv& \frac{\pi}{2}(\f-1)a^2 \mod \pi_1\pi\B.
\end{eqnarray*}
Thus from Sublemma \ref{modsuru}, we obtain
\[
\left(\frac{\f}{2}-1\right) b_x \equiv \Log(f(\pi_n))+\frac{\pi}{2}(\f-1)a^2 \mod \pi_1\pi\B.
\]
\hfill$\square$

From Sublemma \ref{modsuru}, $\theta(b_x)=0$ and there exists an element $b^{\prime}_x \in \Ap$ such that $b_x=b^{\prime}_x v$. By sublemma \ref{factb},
\[
\left\{\left(\frac{\f}{2}-1\right)b_x\right\}\cdot \left(1+\frac{\pi}{2}\right)\equiv \Log(f(\pi_n))\cdot\left(1+\frac{\pi}{2}\right)+\frac{\pi}{2}\left\{(\f-1)a^2\right\}\cdot\left(1+\frac{\pi}{2}\right) \mod \pi_1\pi\B.
\]
Transforming this, we obtain
\[
(\f-v)\left(b_x^{\prime}\cdot\left(1+\frac{\pi}{2}\right)\right) \equiv \Log(f(\pi_n))\cdot\left(1+\frac{\pi}{2}\right)+\pi(\f-1)\left(\frac{a^2}{2}\right) \mod \pi_1\pi\B.
\]
Since the both sides of the above congruence mod $\pi_1\pi\B$ are actually elements in $\Bp$, we have the same congruence mod $\pi_1\pi\Bp=\pi_1\pi\B \cap \Bp$. 
\begin{sublem}\label{liftdekiru}
There exists $c_x \in \Bp$ such that $c_x \equiv b_x^{\prime}\cdot\left(1+\frac{\pi}{2}\right)\mod \pi_1\pi\Bp$ and 
\[
(\f-v)(c_x) = \Log(f(\pi_n))\cdot\left(1+\frac{\pi}{2}\right)+\pi(\f-1)\left(\frac{a^2}{2}\right)
\]
\end{sublem}

(Proof of sublemma \ref{liftdekiru}.)

We show that for any $y \in \pi_1\pi\Bp$, there exists $z$ such that $(\f-v)(z)=y$. For this, it suffices to show the following convergence for any $\pi_1\pi x \ (x \in \Bp)$,  
\[
\left(\frac{\f}{v}\right)^m\left(\frac{\pi_1\pi x}{v} \right):=\left(\frac{\f}{v}\left(\frac{\f}{v}\cdots \left(\frac{\f}{v}\left(\frac{\pi_1\pi x}{v} \right)\right)\cdots\right)\right)\longrightarrow 0\ \ (\text{as $m \rightarrow \infty$}).
\]
In fact, for any $y \in \pi_1\pi\Bp$, a power series $-\sum^{\infty}_{m=0} \left(\frac{\f}{v}\right)^m\left(\frac{y}{v} \right)$ is a solution $z$ of the equation $(\f-v)(z)=y$. If $m=1$, we see that 
\[
\left(\frac{\f}{v}\right)\left(\frac{\pi_1\pi x}{v} \right)= \left(\frac{\f}{v}\right)\left(\pi^2_1 x\right)=\frac{\pi^2 \f(x)}{v}=\pi_1\pi \f(x).
\]
If $m=2$, 
\[
\left(\frac{\f}{v}\right)^2\left(\frac{\pi_1\pi x}{v} \right)=\left(\frac{\f}{v}\right)(\pi_1\pi\f(x))=\pi_1\f(\pi)\f^2(x).
\]
Thus inductively, we have $\left(\frac{\f}{v}\right)^m\left(\frac{\pi_1\pi x}{v} \right)=\pi_1\f^{m-1}(\pi)\f^m(x)$ and $\f^{m-1}(\pi)$ goes to $0$ when $m \rightarrow \infty$ in $\Bp$. Hence we obtain the desired convergence.\hfill$\square$

Dividing the both side of the equation in Sublemma \ref{liftdekiru} by $\pi$, we have
\begin{eqnarray}\label{factc}
(\f-1)\left(\frac{c_x}{\pi_1}-\frac{a^2}{2}\right) = \Log(f(\pi_n))\cdot \left(\frac{1}{\pi}+\frac{1}{2}\right).
\end{eqnarray}
On the other hand, $g \in G_{K_{\cyc}}$ acts on $c_x$ as 
\begin{eqnarray*}
g(c_x v ) \equiv g\left(b_x\cdot\left(1+\frac{\pi}{2}\right)\right)&\equiv&(b_x-\kappa(x)(g) \pi )\cdot\left(1+\frac{\pi}{2}\right)\\
&\equiv& c_xv -\kappa(x)(g)\pi \mod \pi_1\pi\Bp.
\end{eqnarray*}
This implies 
\[
g\left(\frac{c_x}{\pi_1}\right)-\frac{c_x}{\pi_1} \equiv -\kappa(x)(g) \mod \pi_1\Bp.
\]
Since we know $g(a) \equiv a \mod \pi_1\Bp$ from the proof of Sublemma \ref{galactb}, we have 
\begin{eqnarray*}\label{modact}
(g-1)\left(\frac{c_x}{\pi_1}-\frac{a^2}{2}\right)\equiv -\kappa(x)(g) \mod \pi_1\Bp.
\end{eqnarray*}
The above congruence actually yields an equality. In fact, the right hand side $-\kappa(x)(g) \in \Q_2$. For the left hand side, 
\begin{eqnarray*}
(\f-1)\left((g-1)\left(\frac{c_x}{\pi_1}-\frac{a^2}{2}\right)\right)&=&(g-1)\left((\f-1)\left(\frac{c_x}{\pi_1}-\frac{a^2}{2}\right)\right)\\
&=&(g-1)\left(\Log(f(\pi_n))\cdot \left(\frac{1}{\pi}+\frac{1}{2}\right) \right) =0.
\end{eqnarray*}
Here in the second equality, we use (\ref{factc}). There is an exact sequence 
\[
0 \rightarrow \Q_2 \rightarrow \Bc \xrightarrow{\f-1}  \Bc \rightarrow 0.
\] 
Then we see that $(g-1)\left(\frac{c_x}{\pi_1}-\frac{a^2}{2}\right) \in \Q_2$. Since $\Q_2 \cap \pi_1\Bp=0$, we obtain a equality
\[
(g-1)\left(\frac{c_x}{\pi_1}-\frac{a^2}{2}\right)= -\kappa(x)(g).
\]
We now check this $\frac{c_x}{\pi_1}-\frac{a^2}{2}$ is an element in $B_n$. There is an diagram of exact sequences
\[
  \xymatrix{
  0 \ar[r]  & \Q_2 \ar[d]^{id}\ar[r]& B_n \ar[d]^{\mathrm{incl}}\ar[r]^{\f-1}& B_n \ar[d]^{\mathrm{incl}}\ar[r]& 0 \\
  0 \ar[r]& \Q_2 \ar[r] & \Bc \ar[r]^{\f-1}& \Bc \ar[r]& 0. 
  }
\]
Since we have $(\f-1)\left(\frac{c_x}{\pi_1}-\frac{a^2}{2}\right) \in \Bk\subset B_n$ from (\ref{factc}), we can see that $\frac{c_x}{\pi_1}-\frac{a^2}{2} \in B_n$ from the above diagram. Hence, this element $\frac{c_x}{\pi_1}-\frac{a^2}{2}$ is nothing other than the element $\xi_{\kappa(x)} \in B_n$ in Proposition \ref{isocomp}. From (\ref{factc}), we have finished the computation of the first component of $h^1_{\Q_2} \circ \kappa(x)$. We finally compute its second component which we call $\lambda_x(\pi_n)\otimes\e$. Due to Proposition \ref{isocomp}, 
\[
\lambda_x(\pi_n)\otimes\e=(\gh-1)(\xi_x\otimes\e)+\kappa(\gh)\otimes\e.
\]
 We see that 
\begin{eqnarray*}
(\gh-1)(\xi_x\otimes\e)+\kappa(\gh)\otimes\e&=&(\gh-1)\left(\left(\frac{c_x}{\pi_1}-\frac{a^2}{2}\right)\otimes\e\right)+\kappa(\gh)\otimes\e\\
&=& \left(\ch(\g)\gh\left(\frac{c_x}{\pi_1}-\frac{a^2}{2}\right)-\left(\frac{c_x}{\pi_1}-\frac{a^2}{2}\right)\right)\otimes\e+\kappa(\gh)\otimes\e.
\end{eqnarray*}
From Sublemma \ref{liftdekiru}, we have a congruence 
\[
\frac{c_x}{\pi_1} \equiv b_x \cdot \left(\frac{1}{\pi}+\frac{1}{2}\right) \mod \pi\Bp.
\]
Then Sublemma \ref{modsuru} implies
\begin{eqnarray*}
& &\frac{c_x}{\pi_1} \equiv \left(\log a_x +\frac{\pi}{2}a^2\right)\cdot \left(\frac{1}{\pi}+\frac{1}{2}\right)\mod \pi\B\\
& \iff& \frac{c_x}{\pi_1}-\frac{a^2}{2} \equiv \log a_x \cdot \left(\frac{1}{\pi}+\frac{1}{2}\right) \mod \pi\B.
\end{eqnarray*}
On the factor $\frac{1}{\pi}+\frac{1}{2}$, $\gh$ acts as 
\begin{eqnarray}\label{1/pi}
\gh\left(\frac{1}{\pi}+\frac{1}{2}\right)&=&\frac{1}{(1+\pi)^{\ch(\g)}-1}+\frac{1}{2} \\
&\equiv& \frac{1}{\ch(\g)\pi}\cdot\left(1-\frac{\ch(\g)-1}{2}\pi\right)+\frac{1}{2} \mod \pi\B \nonumber \\
&=& \frac{1}{\ch(\g)}\left(\frac{1}{\pi}+\frac{1}{2}\right).\nonumber
\end{eqnarray}
Thus we have 
\begin{eqnarray*}
& & \ch(\g)\gh\left(\frac{c_x}{\pi_1}-\frac{a^2}{2}\right)-\left(\frac{c_x}{\pi_1}-\frac{a^2}{2}\right)\\
&\equiv&  \ch(\g)\gh\left(\log a_x \cdot \left(\frac{1}{\pi}+\frac{1}{2}\right)\right)-\left(\log a_x \cdot \left(\frac{1}{\pi}+\frac{1}{2}\right) \right) \mod \pi\B\\
&\equiv& \ch(\g)\log\gh(a_x)\cdot \frac{1}{\ch(\g)}\left(\frac{1}{\pi}+\frac{1}{2}\right) - \left(\log a_x \cdot \left(\frac{1}{\pi}+\frac{1}{2}\right) \right) \mod \pi\B\\
&=& \left(\log \gh(a_x) -\log a_x\right)\cdot\left(\frac{1}{\pi}+\frac{1}{2}\right).
\end{eqnarray*}
Here, we see that 
\[
\gh(a_x)=\gh\left(\frac{f(\pi_n)}{[\omega_x]}\right)=\frac{\gh(f(\pi_n))}{[\omega_x][\e]^{\kappa(x)(\gh)}}.
\]
Hence, 
\[
\log \gh(a_x) -\log a_x \equiv\gh(\log f(\pi_n))-\log f(\pi_n)-\kappa(x)(\gh)\pi \mod \pi^2\B.
\]
Due to \cite[Lemma 2.2.1]{Ben}, 
\[
\gh(\log f(\pi_n))-\log f(\pi_n) \equiv \frac{\ch(\g)-1}{2^n}D\log f(\pi_n)\cdot\pi \mod \pi^2\Bp.
\]
This implies a congruence mod $\pi\Bp$
\begin{eqnarray*}
& &\ch(\g)\gh\left(\frac{c_x}{\pi_1}-\frac{a^2}{2}\right)-\left(\frac{c_x}{\pi_1}-\frac{a^2}{2}\right) \\
&\equiv& \left(\gh(\log f(\pi_n))-\log f(\pi_n)-\kappa(x)(\gh)\pi\right)\cdot\left(\frac{1}{\pi}+\frac{1}{2}\right)\ \mod \pi\Bp,
\end{eqnarray*}
where we use $\pi\B\cap\Bc=\pi\Bp$. Thus we obtain
\begin{eqnarray*}
\lambda_x(\pi_n)&=&\left(\ch(\g)\gh\left(\frac{c_x}{\pi_1}-\frac{a^2}{2}\right)-\left(\frac{c_x}{\pi_1}-\frac{a^2}{2}\right)\right)+\kappa(\gh)\\
&\equiv& \left(\frac{\ch(\g)-1}{2^n}D\log f(\pi_n)\cdot\pi-\kappa(x)(\gh)\pi\right)\cdot\left(\frac{1}{\pi}+\frac{1}{2}\right)+\kappa(\gh)\mod \pi\Bp\\
&\equiv&  \frac{\ch(\g)-1}{2^n}D\log(f(\pi_n)) \mod \pi\Bp.
\end{eqnarray*}
However, since $\lambda_x(\pi_n) \in \Bk$ and $\pi\Bp \cap \Bk=\pi\Ok{[[\pi_n]]}\tens{\mathbb{Z}_2}\mathbb{Q}_2$, the above congruence mod $\pi\Bp$ is in fact the one mod $\pi\Ok{[[\pi_n]]}\tens{\mathbb{Z}_2}\mathbb{Q}_2$ and hence $\lambda_x(\pi_n) \in \Ok{[[\pi_n]]}\tens{\mathbb{Z}_2}\mathbb{Q}_2$. Thus we finally obtain the claim of Lemma \ref{tensorQ}.\hfill$\square$ 
 \begin{lem}\label{slide}
 There exists a power series $Y_x(\pi_n) \in \frac{1}{2}\Ak$ such that 
 \[
 (\f-1)Y_x(\pi_n)=\frac{1}{2}\Log(f(\pi_n)) \
 \]
 \end{lem}
(Proof of Lemma \ref{slide})
\\ 
 Since $x \in U_{K_n}^1$, we have $f(\pi_n) \in 1+\pi_n\Ok{[[\pi_n]]}$ and $\Log(f(\pi_n))=\left(\frac{\f}{p}-1\right)\log f(\pi_n) \in \pi_n\Ok{[[\pi_n]]}$. We define
 \[
 Y_x(\pi_n):=-\sum^{\infty}_{i=0}\f^i\left(\Log(f(\pi_n))\cdot \frac{1}{2}\right).
 \]
Note that this $Y_x(\pi_n)$ is a well-defined element in $\frac{1}{2}\Ak$ since $\f^i(\pi_n) \rightarrow 0$ as $i \rightarrow \infty$ in $\Bk$. We can see that $Y_x(\pi_n)$ satisfies $(\f-1)(Y_x(\pi_n))=\frac{1}{2}\Log(f(\pi_n))$\hfill$\square$\\

From Lemma \ref{slide}, we have a 1-coboundary of the complex $C^{\bullet}(\Bk(1))$
\[
\left[(\f-1)(Y_x(\pi_n)\otimes\e), (\g-1)(Y_x(\pi_n)\otimes\e)  \right] = \left[\Log(f(\pi_n))\cdot \frac{1}{2} \otimes \e, (\ch(\g)-1)Y_x(\pi_n)\otimes\e\right]
\]
Subtracting this 1-coboundary from the result in Lemma \ref{tensorQ}, we obtain
\begin{eqnarray}\label{kettensorQ}
\hspace{10mm}h^1_{\Q_2}\circ\kappa(x)= \left[\Log(f(\pi_n))\cdot\frac{1}{\pi}\otimes \e,\ \lambda_x (\pi_n) \otimes \e -(\ch(\g)-1)Y_x(\pi_n)\otimes\e\right].
\end{eqnarray}
Thus the first component of the above representative for $h^1_{\Q_2}\circ\kappa(x)$ is actually an element in $\Ak(1)$. We show that so does the second component.
\begin{lem}\label{nodenom}
We have $\lambda_x (\pi_n) \in \Ok{[[\pi_n]]}$, hence
\[
\lambda_x(\pi_n) \equiv \frac{\ch(\g)-1}{2^n}D\log f(\pi_n)\ \ \mathrm{mod}\ \pi\Ok{[[\pi_n]]}.
\]
\end{lem}
(Proof of Lemma \ref{nodenom})

From (\ref{kettensorQ}), the 1-cocycle condition says
\begin{eqnarray*}
& &(\g-1)\left(\Log(f(\pi_n))\cdot\frac{1}{\pi}\otimes \e\right)=(\f-1)\left(\lambda_x (\pi_n) \otimes \e -(\ch(\g)-1)Y_x(\pi_n)\otimes\e\right)\\
&\hspace{-10mm}\iff & (\f-1)(\lambda_x (\pi_n))=\ch(\g)\g\left(\Log(f)\cdot\frac{1}{\pi}\right)-\Log(f)\cdot\frac{1}{\pi}+(\f-1)(\ch(\g)-1)Y_x.
\end{eqnarray*}
Then we can see that $(\f-1)(\lambda_x (\pi_n)) \in \Ak$. While, there is a commutative diagram
\[
  \xymatrix{
  0 \ar[r]  & \Z_2 \ar[d]^{\mathrm{incl}}\ar[r]& A_n \ar[d]^{\mathrm{incl}}\ar[r]^{\f-1}& A_n \ar[d]^{\mathrm{incl}}\ar[r]& 0 \\
  0 \ar[r]& \Q_2 \ar[r] & B_n \ar[r]^{\f-1}& B_n\ar[r]& 0 
  }
\]
which implies that there exists $r \in \Q_2$ such that $\lambda_x (\pi_n)-r \in A_n$. However, from Lemma \ref{tensorQ}, we have 
\[
\lambda_x(\pi_n) \equiv \frac{\chi(\g)-1}{2^n}D\log f(\pi_n)\ \ \mathrm{mod}\ \pi\Ok{[[\pi_n]]}\tens{\Z_2}\Q_2.
\]
In other words, we can see that $\lambda_x(\pi_n) = (\text{an element in $\Ak$})+(\text{terms divisible by $\pi$})$. Hence $r$ must be $0$ and $\lambda_x(\pi_n) \in A_n \cap (\Ok{[[\pi_n]]}\tens{\Z_2}\Q_2)=\Ok{[[\pi_n]]}$.\hfill$\square$

We finally prove Proposition \ref{key}.

(Proof of Proposition \ref{key})

There is a commutative diagram
\[
  \xymatrix{
 (U_{K_n}^1)^f \ar[r]^{\hspace{-7mm}\kappa} & (H^1(K_n, \Z_2(1)))^f \ar[d]^{\iota}\ar[r]^{h^1}_{\sim}  & (H^1_{\Phi\G}(\Ak(1)))^f \ar[d]^{\iota_{\Phi\G}} \\
  &H^1(K_n, \Q_2(1)) \ar[r]^{h^1_{\Q_2}}_{\sim} &  H^1_{\Phi\G}(\Bk(1)), 
  }
\]
where $(M)^f$ denotes the torsion-free part of a $\Z_2$-module $M$. Note also that $\iota, \iota_{\Phi\G}$ are the homomorphisms which induced by inclusions. Since we consider only torsion-free parts of $\Z_2$-modules in the first row, the vertical arrows $\iota, \iota_{\Phi\G}$ are injective. From (\ref{kettensorQ}) and Lemma \ref{nodenom}, for any $x \in (U_{K_n}^1)^f$, we have 
\[
h^1_{\Q_2}\circ\kappa(x)= \left[\Log(f(\pi_n))\cdot\frac{1}{\pi}\otimes \e,\ \lambda_x (\pi_n) \otimes \e -(\chi_{\cyc}(\g)-1)Y_x(\pi_n)\otimes\e\right],
\]
and the first and second components of the above representative are in $\Ak(1)$. Thus the pair of elements in $\Ak$
\[
\left(\Log(f(\pi_n))\cdot\frac{1}{\pi}\otimes \e,\ \lambda_x (\pi_n) \otimes \e -(\chi_{\cyc}(\g)-1)Y_x(\pi_n)\otimes\e\right)
\]
 also defines a cohomology class in $H^1_{\Phi\G}(\Ak(1))$ which maps to $h^1_{\Q_2}\circ\kappa(x)$ under $\iota_{\Phi\G}$. By the commutativity of the above diagram and the injectivity of $\iota_{\Phi\G}$, we have
 \[
h^1\circ\ka (x)=\left[\Log(f(\pi_n))\cdot\frac{1}{\pi} \otimes \e,\ \lambda_x(\pi_n) \otimes\e -(\chi_{\cyc}(\g)-1)Y_x(\pi_n)\otimes \e\right].
 \]
This completes the proof of Proposition \ref{key}. \hfill$\square$
\section{Calculation of the Hilbert symbol} In this section, we calculate the Hilbert symbol and give an explicit formula following the strategy we mentioned in Subsection \ref{str}.
\subsection{Computation of the cup product $\cup_{\Phi\G}$}
\begin{lem}\label{cupcomp}
Let $x, y \in (U_{K_n}^1)^f$. There is a power series $H_{x, y} \in \Ak$ such that $(h^1\circ\ka (x)) \cup_{\Phi\G}(h^1\circ\ka (y))=[H_{x, y} \otimes \varepsilon^2]$ and 
\begin{eqnarray*}
H_{x, y}&\equiv&\frac{\ch(\g)-1}{2^n}\left(D\log f\cdot\mathfrak{L}(g)-\mathfrak{L}(f)\f(D\log f)\right)\cdot\frac{1}{\pi} \\
& & +(\ch(\g)-1)\left(\mathfrak{L}(f)\f(Y_y)-Y_x\mathfrak{L}(g)\right)\cdot \frac{1}{\pi}\ \ \mathrm{mod}\ \Ok{[[\pi_n]]}.
\end{eqnarray*}
Here $f(\pi_n), g(\pi_n)\in \Ok{[[\pi_n]]}$ are power series which satisfy $f(\zeta_{2^n}-1)=x, g(\zeta_{2^n}-1)=y$.
\end{lem}
(Proof of Lemma \ref{cupcomp})

Using Proposition \ref{key}, we have 
\begin{eqnarray*}
h^1\circ\ka (x)&=&\left[\Log(f(\pi_n))\cdot\frac{1}{\pi} \otimes \e,\ \lambda_x(\pi_n) \otimes\e -(\ch(\g)-1)Y_x(\pi_n)\otimes \e\right],\\
h^1\circ\ka (y)&=&\left[\Log(g(\pi_n))\cdot\frac{1}{\pi} \otimes \e,\ \lambda_y(\pi_n) \otimes\e -(\ch(\g)-1)Y_y(\pi_n)\otimes \e\right].
\end{eqnarray*}
From Proposition \ref{cup}, we can compute the cup product as $(h^1\circ\ka (x)) \cup_{\Phi\G}(h^1\circ\ka (y))=[H_{x, y} \otimes \varepsilon^2]$, where
\begin{eqnarray*}
H_{x, y}&=&\left(\lambda_x -(\ch(\g)-1)Y_x\right)\cdot\ch(\g)\g\left(\Log(g)\cdot\frac{1}{\pi}\right) \\
& &-\left(\Log(f)\cdot\frac{1}{\pi}\right)\cdot\f\left(\lambda_y -(\ch(\g)-1)Y_y\right).
\end{eqnarray*}
As we saw in (\ref{1/pi}), we have
\[
\g\left(\frac{1}{\pi}\right) \equiv \frac{1}{\ch(\g)\pi}\ \  \mathrm{mod}\ \Ok{[[\pi_n]]}, 
\]
and from \cite[Lemma 2.2.1]{Ben}, for $F(X) \in \Ok{[[X]]}$, we also have
\[
\g(F(\pi_n)) \equiv F(\pi_n) \ \ \mathrm{mod}\ \pi\Ok{[[\pi_n]]}.
\]
These congruences implies
\begin{eqnarray*}
H_{x, y}&\equiv&(\lambda_x-(\ch(\g)-1)Y_x)\mathfrak{L}(g)\cdot\frac{1}{\pi} \\
& &-\Log(f)\cdot\f\left(\lambda_y -(\ch(\g)-1)Y_y\right)\cdot \frac{1}{\pi}\ \ \mathrm{mod}\ \Ok{[[\pi_n]]}.
\end{eqnarray*}
Here, from Proposition \ref{key}, we know
\[
\lambda_x(\pi_n) \equiv \frac{\ch(\g)-1}{2^n}D\log f,\ \ \lambda_y(\pi_n) \equiv \frac{\ch(\g)-1}{2^n}D\log g\ \ \mathrm{mod}\ \pi\Ok{[[\pi_n]]}.
\]
Then we obtain
\begin{eqnarray*}
H_{x, y}&\equiv&\frac{\ch(\g)-1}{2^n}\left(D\log f\cdot\mathfrak{L}(g)-\mathfrak{L}(f)\f(D\log f)\right)\cdot\frac{1}{\pi} \\
& & +(\ch(\g)-1)\left(\mathfrak{L}(f)\f(Y_y)-Y_x\mathfrak{L}(g)\right)\cdot \frac{1}{\pi}\ \ \mathrm{mod}\ \Ok{[[\pi_n]]}.
\end{eqnarray*}
\hfill$\square$
\subsection{Explicit formula for the Hilbert symbol}
We finally compute the image of $(h^1\circ\ka (x)) \cup_{\Phi\G}(h^1\circ\ka (y))$ under $\overline{\mathrm{TR}_{K_n}}$ and complete the calculation of the Hilbert symbol.
\begin{thm}
For $x, y \in U_{K_n}^1$, 
\begin{eqnarray*}\label{mainthm}
& &[x, y]_{K_n}\\ 
&\hspace{-3mm}=&\hspace{-2mm}-(1+2^{n-1})\Tr_{K/\Q_2}\left(\Res_{\pi_n}\left(D\log f\cdot\Log(g)-\Log(f)\f(D\log(g))\right)
\frac{d\pi_n}{\pi(1+\pi_n)}\right)\\
& &-2^n\Tr_{K/\Q_2}\left(\Res_{\pi_n}\left(\Log(f)\f(Y_y)-Y_x\Log(g)\right)\frac{d\pi_n}{\pi(1+\pi_n)}\right).
\end{eqnarray*}
Here power series $f(\pi_n), g(\pi_n)$ are the same as in Lemma \ref{cupcomp}.
\end{thm}

(Proof of Theorem \ref{mainthm})

First we show the theorem for $x, y \in (U_{K_n}^1)^f$. All we have to do is just computing $\mathrm{TR}_n(H_{x, y}\otimes\e)\ \mathrm{mod}\ 2^n$. By the fact that elements in $\Ok{[[\pi_n]]}$ have no residue and Lemma \ref{cupcomp}, we have
\begin{eqnarray*}
& &\mathrm{TR}_n(H_{x, y}\otimes\e)\\ 
&=&\mathrm{TR}_n\left(\frac{\ch(\g)-1}{2^n}\left(D\log f\cdot\mathfrak{L}(g)-\mathfrak{L}(f)\f(D\log f)\right)\cdot\frac{1}{\pi}\right) \\
&& +\mathrm{TR}_n\left((\ch(\g)-1)\left(\mathfrak{L}(f)\f(Y_y)-Y_x\mathfrak{L}(g)\right)\cdot \frac{1}{\pi}\right)\\
&=&\hspace{-2mm}-\frac{\ch(\g)-1}{\log(\ch(\g))}\Tr_{K/\Q_2}\left(\Res_{\pi_n}\left(D\log f\cdot\mathfrak{L}(g)-\mathfrak{L}(f)\f(D\log f)\right)\frac{d\pi_n}{\pi(1+\pi_n)}\right)
\\
& &-2^n\frac{\ch(\g)-1}{\log(\ch(\g))}\Tr_{K/\Q_2}\left(\Res_{\pi_n}\left(\Log(f)\f(Y_y)-Y_x\Log(g)\right)\frac{d\pi_n}{\pi(1+\pi_n)}\right).
\end{eqnarray*}
On the other hand, we can see that
\[
\frac{\ch(\g)-1}{\log(\ch(\g))}\equiv 1+\frac{1}{2}(\ch(\g)-1) \pmod{2^n}.
\]
Since $\g$ is a topological generator of the Galois group $\Gamma_n$, there exists $u \in \Z_2^{\times}$ such that $\ch(\g)-1=2^n u$. Then we have $\frac{1}{2}(\ch(\g)-1)=2^{n-1}u\equiv 2^{n-1}\pmod{2^n}$ because $\Z_2^{\times}=\langle-1, 5\rangle$. Thus we obtain Theorem \ref{mainthm} when $x, y \in (U_{K_n}^1)^f$.

Next we consider the case that one of $x$ and $y$ is not in $(U_{K_n}^1)^f$. Since $U_{K_n}^1=\langle \zeta_{2^n} \rangle \oplus (U_{K_n}^1)^f$, it suffices to consider the case when $y=\zeta_{2^n}$. in the following, we use the Artin-Hasse formula and some facts from \cite{Ben} on power series.
\begin{thm*}[Artin-Hasse, \cite{AH}]\label{AH}
For $y \in U_{\Q_2(\zeta_{2^n})}^1$, 
\[
[x, \zeta_{p^n}]_{\Q_2(\zeta_{2^n})}=-\frac{1+2^{n-1}}{2^n}\Tr_{\Q_2(\zeta_{2^n})/\Q_2}(\log x).
\]
\end{thm*}
\begin{lem}[Proposition 2.2.1, \cite{Ben}]\label{fact1}
For any $F(X) \in \Ok{[[X]]}$, 
\[
\Res_{\pi_n}\left(F(\pi_n)\frac{d\pi_n}{\pi(1+\pi_n)}\right)=\frac{1}{2^n}\sum_{\zeta\in\mu_{2^n}}F(\zeta-1).
\]
\end{lem}
\begin{lem}[Lemma 2.2.5.1, \cite{Ben}]\label{fact2}
Let $x \in U_{K_n}$ and $f(X)\in\Ok{[[X]]}$ which satisfies $f(\zeta_{2^n}
-1)=x$. Then
\[
\Tr_{K_n/\Q_2}\log x = -\Tr_{K/\Q_2}\left(\sum_{\zeta\in\mu_{2^n}}\mathfrak{L}(f)(\zeta-1)\right).
\]
\end{lem}
We verify the validity of Theorem \ref{mainthm} for $x\in U_{K_n}^1$ and $y=\zeta_{2^n}$. First we compute the Hilbert symbol via the Artin-Hasse formula. We see that
\begin{eqnarray*}
(x, \zeta_{2^n})_{K_n}=\frac{\rho_{K_n}(x)(\zeta_{2^{2n}})}{\zeta_{2^{2n}}}=\frac{\rho_{K_n}(x)\mid_{\Q_2(\zeta_{2^n})^{\mathrm{ab}}}(\zeta_{2^{2n}})}{\zeta_{2^{2n}}}&=&\frac{\rho_{\Q_2(\zeta_{2^n})}(\mathrm{N}_{K_n/\Q_2(\zeta_{2^n})}(x))(\zeta_{2^{2n}})}{\zeta_{2^{2n}}}\\
&=& (\zeta_{2^n}, \mathrm{N}_{K_n/\Q_2(\zeta_{2^n})}(x))_{\Q_2(\zeta_{2^n})},
\end{eqnarray*}
where $\mathrm{N}_{K_n/\Q_2(\zeta_{2^n})}$ denotes the field norm of the extension $K_n/\Q_2(\zeta_{2^n})$. Then the Artin-Hasse formula implies
\begin{eqnarray*}
[x, \zeta_{2^n}]_{K_n}=[\mathrm{N}_{K_n/\Q_2(\zeta_{2^n})}(x), \zeta_{2^n}]_{\Q_2(\zeta_{2^n})}&=&-\frac{1+2^{n-1}}{2^n}\Tr_{\Q_2(\zeta_{2^n})/\Q_2}(\log(\mathrm{N}_{K_n/\Q_2(\zeta_{2^n})}(x)))\\
&=&-\frac{1+2^{n-1}}{2^n}\Tr_{K_n/\Q_2}(\log x).
\end{eqnarray*}
Next we compute the right-hand side of the formula in Theorem \ref{mainthm}. When $y=\zeta_{2^n}$, we can take $g(\pi_n)=\pi_n-1$ to get $\mathfrak{L}(g(\pi_n))=\left(\frac{\f}{2}-1\right)\log(1+\pi_n)$. Hence by Lemma \ref{fact1} and the definition of the power series $Y_y(\pi_n)$, we have 
\begin{eqnarray*}
& &\Res_{\pi_n}\left(\Log(f)\f(Y_y)-Y_x\Log(g)\right)\frac{d\pi_n}{\pi(1+\pi_n)}\\
& = & \frac{1}{2^n}\sum_{\zeta\in\mu_{2^n}}\left(\Log(f(X))\f(Y_y(X))-Y_x(X)\Log(g(X))\right)\mid_{X=\zeta-1}=0.\\
\end{eqnarray*}
Similarly, we also have $\Res_{\pi_n}\left(D\log f\cdot\Log(g)\right)=0$. Thus we can see that 
\begin{eqnarray*}
&\hspace{-2mm}&-(1+2^{n-1})\Tr_{K/\Q_2}\left(\Res_{\pi_n}\left(D\log f\cdot\Log(g)-\Log(f)\f(D\log(g))\right)
\frac{d\pi_n}{\pi(1+\pi_n)}\right)\\
& &-2^n\Tr_{K/\Q_2}\left(\Res_{\pi_n}\left(\Log(f)\f(Y_y)-Y_x\Log(g)\right)\frac{d\pi_n}{\pi(1+\pi_n)}\right)\\
&=&(1+2^{n-1})\Tr_{K/\Q_2}\left(\Res_{\pi_n}\left(\Log(f)\f(D\log(g))\right)
\frac{d\pi_n}{\pi(1+\pi_n)}\right)\\
& = & (1+2^{n-1})\Tr_{K/\Q_2}\left(\frac{1}{2^n}\sum_{\zeta\in\mu_{2^n}}\left(\Log(f)\f(D\log(g))\right)\mid_{X=\zeta-1}\right)\\
& = & \frac{1+2^{n-1}}{2^n}\Tr_{K/\Q_2}\sum_{\zeta\in\mu_{2^n}}\left(\Log(f(\zeta-1))\right)\\
&=&-\frac{1+2^{n-1}}{2^n}\Tr_{K/\Q_2}\left(\Tr_{K_n/\Q_2}\log x\right)=-\frac{1+2^{n-1}}{2^n}\Tr_{K_n/\Q_2}\left(\log x\right)
\end{eqnarray*}
Here we use  
\[
\f(D\log(g(X)))=\f((1+X)\frac{d}{dX}\log (1+X))=1
\]
in the third equality and Lemma \ref{fact2} in the fourth equality. This completes the proof.

\hfill$\square$


\begin{thebibliography}{99}
\bibitem{Abs} V. A. Abrashkin, The field of norms functor and the Br\"uckner-Vostokov formula, Math. Ann. 308 (1997), 5-19.
\bibitem{Ai} K. Aichi, On Iwasawa's explicit formula for the norm residue symbol, Memoirs of The Faculty of Science, Kyushu University. Series A, Mathematics 26 (1972), 139-148.
\bibitem{AH}  E. Artin, H. Hasse, Die beiden Erg\"anzungss\"atze zum Reziprozit\"atsgesetz der $l^n$-ten Potenzreste im K\"orper der $l^n$-ten Einheitswurzeln, Abh. Math. Semin. Univ. Hamb. 6 (1928) 146-162.
\bibitem{Ben} D.Benois, On Iwasawa theory of crystalline representations, Duke Math. J. 104 (2000), no. 2, 211-267.
\bibitem{Col1} R. Coleman, Division values in local fields, Invent. Math. 53 (1979), 91-116
\bibitem{Col2} R. Coleman,The dilogarithm and the norm residue symbol, Bull. Soc. Math. France 109 (1981), 373-402.
\bibitem{dS} E. de Shalit, The explicit reciprocity law in local class field theory, Duke Mathematical Journal 53 (1986), 163-176.
\bibitem{Fl} J. Fl\'orez, Explicit reciprocity laws for higher local fields, Journal of Number Theory. 213 (2020), 400-444.
\bibitem{F} Fontaine, J.-M.,“Repr\'esentations $p$-adiques des corps locaux, I”in The Grothendieck Festschrift, Vol. 2, Progr. Math. 87, Birkh\"auser, Boston, 1990, 249-309.
\bibitem{FO}  J.-M. Fontaine,Yi Ouyang, Theory of $p$-adic Galois-representations: staff.ustc.edu. cn/ yiouyang/galoisrep.pdf
\bibitem{FV} I.B.Fesenko and S.V.Vostokov, Local Fields and Their Extensions, Amer. Math. Soc., Providence, RI, 2002.
\bibitem{H1}  Laurent Herr, Sur la cohomologie galoisienne des corps $p$-adiques, Bull. Soc. Math. France 126 (1998), no. 4, 563-600. 
\bibitem{H2} Laurent Herr, Une approche nouvelle de la dualit\'e locale de Tate, Math. Ann. 320 (2001), no. 2, 307-337. 
\bibitem{Iwa1} K. Iwasawa, 1986. Local Class Field Theory. Oxford Science Publications.
\bibitem{Iwa2} K. Iwasawa, On Explicit Formulas for the Norm Residue Symbol, J. Math. Soc. Japan 20 (1968), 151-165.
\bibitem{Ka} K. Kato, The explicit reciprocity law and the cohomology of Fontaine-Messing, Bull. Soc. Math. France 119 (1991), 397-441.
\bibitem{Ku} M. Kurihara, The exponential homomorphism for the Milnor K-groups and an explicit reci-
procity law, J. Reine Angew. 498 (1998) 201-221.
\bibitem{NSW}J.Neukirch, A.Schmidt, K.Wingberg, Cohomology of Number Fields,
Springer-Verlag (2008).
\bibitem{W}  A. Wiles, Higher explicit reciprocity laws, Ann. Math. 107 (2) (1978) 235-254.
\bibitem{V} Sergei V. Vostokov, Explicit formulas for the Hilbert symbol, Invitation to higher local fields (M\"unster, 1999), Geom. Topol. Monogr., vol. 3, Geom. Topol. Publ., Coventry, 2000, pp. 81-89. 
\bibitem{Z}  A. N. Zinoviev, Generalized Artin-Hasse and Iwasawa formulas for the Hilbert symbol in a
high-dimensional local field, JMS 124 (1) (2004) 4806-4821.
\end{thebibliography}
\end{document}